\documentclass{amsart}
\usepackage{amsfonts}
\usepackage{amsmath}
\usepackage{amssymb}
\usepackage{amsthm}

\newtheorem{theorem}{Theorem}
\newtheorem{prop}{Proposition}
\newtheorem{lemma}{Lemma}
\newtheorem{rem}{Remark}

\begin{document}
\author{Mark Pankov}
\title[Order preserving transformations of the Hilbert grassmannian]
{Order preserving transformations of the Hilbert grassmannian:
complex case}
\address{Department of Mathematics and Information Technology,
University of Warmia and Mazury, {\. Z}olnierska 14A, 10-561
Olsztyn, Poland} \email{pankov@matman.uwm.edu.pl}

\begin{abstract}
Let $H$ be a separable complex Hilbert space.
Denote by ${\mathcal G}_{\infty}(H)$ the Grassmannian consisting
of closed linear subspaces with infinite dimension and codimension.
This Grassmannian is partially ordered by the inclusion relation.
We show that every continuous order preserving bijective transformation of ${\mathcal G}_{\infty}(H)$
is induced by an invertible bounded semi-linear operator.
\end{abstract}

\maketitle

\section{Introduction}
Let $H$ be a separable (real or complex) Hilbert space.
We write ${\mathcal G}(H)$ for the lattice of closed linear subspaces of $H$.
The group of invertible bounded linear operators
acts on ${\mathcal G}(H)$ and the orbits of this action is called
{\it Grassmannians}:
\begin{enumerate}
\item[$\bullet$]
${\mathcal G}_{k}(H)$ consists of $k$-dimensional linear subspaces,
\item[$\bullet$]
${\mathcal G}^{k}(H)$ consists of closed linear subspaces with
codimension $k$,
\item[$\bullet$]
${\mathcal G}_{\infty}(H)$ consists of closed linear subspaces
with infinite dimension and codimension.
\end{enumerate}
Every closed linear subspace $S\subset H$ can be identified with
the orthogonal (self-adjoint) projection $p_{S}:H\to S$,
then ${\mathcal G}(H)$ and all our Grassmannian are closed subsets in
the Banach algebra with the standard operator norm.
In other words, these
are complete metric spaces with the distance
$||P_{S}-P_{U}||$
(we refer \cite{Shubin} for some elementary topological properties of the Grassmannians).

Let us consider the Grassmannian ${\mathcal G}_{\infty}(H)$
which is partially ordered by the inclusion relation.

In the real case,
Mackey's result \cite{Mackey}
states that every automorphism of the lattice ${\mathcal G}(H)$
is induced by an invertible bounded linear operator,
and it was shown by author \cite{Pankov} that
every order preserving bijective transformation
of ${\mathcal G}_{\infty}(H)$
can be extended to an automorphism of the lattice.
In the present paper complex versions of these results will be given.

Now suppose that our Hilbert space is complex.
We say that $A:H\to H$ is a {\it semi-linear operator}
if
$$A(x+y)=A(x)+A(y)$$
for all $x,y\in H$
and there exists an automorphism $\sigma: \mathbb{C}\to  \mathbb{C}$
such that
$$A(ax)=\sigma(a)A(x)$$
for all $a\in \mathbb{C}$ and $x\in H$.
If $A$ is a non-zero bounded semi-linear operator then the associated automorphism
$\sigma: \mathbb{C}\to  \mathbb{C}$ is continuous.
Since $\sigma(q)=q$ for every $q\in \mathbb{Q}$
and $\sigma(i)=\pm i$,
one of the following possibilities is realized:
\begin{enumerate}
\item[$\bullet$]
$\sigma$ is identical and $A$ is linear,
\item[$\bullet$]
$\sigma$ is the complex conjugate mapping and $A(ax)=\bar{a}A(x)$
for all $a\in \mathbb{C}$, $x\in H$.
\end{enumerate}

Note that non-continuous automorphisms of $\mathbb{C}$ exist.

Every invertible bounded semi-linear operator induces a
continuous order preserving bijective transformation of
${\mathcal G}_{\infty}(H)$.
Conversely, the following statement will be proved.

\begin{theorem}
Every continuous order preserving bijective transformation of
${\mathcal G}_{\infty}(H)$ is induced by
an invertible bounded semi-linear operator.
\end{theorem}

In the real case \cite{Pankov}, we do not require continuity.

\section{Hilbert projective spaces}
Let $H$ be, as in the previous section, a separable (real or complex)
Hilbert space.
Consider the associated {\it Hilbert projective space}
$\Pi_{H}$ whose points are $1$-dimensional linear subspaces and whose lines are defined by
$2$-dimensional linear subspaces of $H$.
It is clear that lines are closed and all closed subspaces of $\Pi_{H}$
are induced by closed linear subspaces of $H$.
The projective space whose points are closed hyperplanes
and whose lines are defined by closed linear subspaces of codimension $2$
will be called the {\it dual Hilbert projective space}
and denoted by $\Pi^{*}_{H}$.
The mapping $S\to S^{\perp}$ ($S^{\perp}$ is the orthogonal complement)
induces an isometry between ${\mathcal G}_{k}(H)$ and ${\mathcal G}^{k}(H)$.
By this mapping,
$\Pi_{H}$ and $\Pi^{*}_{H}$ are canonically isomorphic
as metric projective spaces.

It was proved in  \cite{Mackey} that
every collineation of $\Pi_{H}$ preserving the class of closed hyperplanes
("closed" collineation)
is induced by an invertible bounded linear operator
if $H$ is real.
If our Hilbert space is complex then
the following "weak" version of Mackey's result is true.

\begin{theorem}
Every continuous collineation of $\Pi_{H}$
is induced by an invertible bounded semi-linear operator.
\end{theorem}

To prove Theorem 2 we use the following lemma.

\begin{lemma}
Let $A:H\to H$ be an invertible semi-linear operator
such that the associated automorphism of $\mathbb{C}$
is identical or the complex conjugate mapping.
If $A$ preserves ${\mathcal G}^{1}(H)$ then it is bounded.
\end{lemma}

\begin{proof}
Similar to the proof of Lemma B in \cite{Mackey}.
\end{proof}

\begin{proof}[Proof of Theorem 2]
Let $f:{\mathcal G}_{1}(H)\to {\mathcal G}_{1}(H)$
be a continuous collineation of $\Pi_{H}$.
By the Fundamental Theorem of Projective Geometry \cite{Baer},
there exists an invertible semi-linear operator
$A:H\to H$ such that
$$f(P)=A(P)$$
for all $P\in {\mathcal G}_{1}(H)$.
For the associated automorphism $\sigma:{\mathbb C}\to {\mathbb C}$
one of the following possibilities is realized:
\begin{enumerate}
\item[(1)]
the restriction of $\sigma$ to ${\mathbb Q} + {\mathbb Q}i$
is identical,
\item[(2)]
$\sigma$ transfers $p+qi$ to $p-qi$ for each $p,q\in {\mathbb Q}$.
\end{enumerate}
Suppose that $x,y\in H$ are linearly independent
and $g$ is the continuous collineation of $\Pi_{H}$ induced by
an invertible bounded semi-linear operator transferring
$A(x)$ and $A(y)$ to $x$ and $y$, respectively;
we require that the associated automorphism $\varphi$
is identical or the complex conjugate mapping in the case (1) or (2), respectively;
in other words,
the restriction of $\varphi\sigma$ to ${\mathbb Q} + {\mathbb Q}i$
is identical.
Then $gf$ preserves the line of $\Pi_{H}$ joining $\langle x \rangle$ and
$\langle y \rangle$, and it transfers $\langle x+ay \rangle$
to $\langle x+\varphi\sigma(a)y \rangle$.
Since
$$\langle x+by \rangle,\;\;b\in{\mathbb Q} + {\mathbb Q}i$$
form an everywhere dense subset in the line,
the restriction of $gf$ to the line is identical. Hence $\sigma=\varphi$.
We apply Lemma 1 and finish the proof.
\end{proof}

\begin{rem}{\rm
It follows from Theorem 2 that every continuous automorphism
of the lattice ${\mathcal G}(H)$ is induced by an invertible
bounded semi-linear operator.
}\end{rem}

\section{Proof of Theorem 1}
We will use the following trivial remark concerning the automorphisms of ${\mathcal G}(H)$
induced by invertible bounded semi-linear operators.
If $A:H\to H$ is an invertible bounded semi-linear operator
then the transformation
$$S\to (A(S^{\perp}))^{\perp}$$
is induced by the operator $(A^{*})^{-1}$,
where $A^{*}$ is the adjoint operator.
Recall that $A^{*}:H\to H$ is defined by the formula
$$(A^{*}x, y)=(x,Ay)\;\;\;\;\;\forall\;x,y\in H$$
or
$$(A^{*}x, y)=\overline{(x,Ay)}\;\;\;\;\;\forall\;x,y\in H$$
if the associated automorphism of $\mathbb{C}$ is identical or
the complex conjugate mapping, respectively.

Let $f:{\mathcal G}_{\infty}(H)\to{\mathcal G}_{\infty}(H)$
be a continuous order preserving bijection.

\begin{lemma}
Let $S\in {\mathcal G}_{\infty}(H)$
and ${\mathcal X}$ be the set of all elements of ${\mathcal G}_{\infty}(H)$
contained in $S$.
There exists an invertible bounded semi-linear operator $A:S\to f(S)$ such that
$$f(U)=A(U)\;\;\;\;\;\forall\;U\in {\mathcal X}.$$
\end{lemma}

\begin{proof}
We suppose that $f(S)=S$,
since in the general case we can take an invertible bounded
linear operator $C:H\to H$
which sends $f(S)$ to $S$ and consider the transformation $U\to C(f(U))$.
The restriction of $f$ to ${\mathcal G}^{1}(S)$
is a continuous collineation of $\Pi^{*}_{S}$.
Then
$$P\to (f(P^{\perp}\cap S))^{\perp}\cap S\;\;\;\;\;\forall\;P\in {\mathcal G}_{1}(S)$$
is a continuous collineation of $\Pi_{S}$.
By Theorem 2, it is induced by an invertible bounded semi-linear operator
$B:S\to S$.
Hence the restriction of $f$ to ${\mathcal G}^{1}(S)$
is induced by the operator $$A:=(B^{*})^{-1}.$$
Let $[U]^{1}$ be the set of all elements of ${\mathcal G}^{1}(S)$
containing $U\in {\mathcal X}$.
Since $f$ is order preserving,
$$f([U]^{1})=[f(U)]^{1};$$
on the other hand,
$$f([U]^{1})=A([U]^{1})=[A(U)]^{1}.$$
Thus $f(U)=A(U)$ and $A$ is as required.
\end{proof}

\begin{rem}{\rm
In the real case,
an order preserving transformation is not assumed to be continuous
and the proof of the analogous lemma (Lemma 3 in \cite{Pankov})
is more complicated.
}\end{rem}

\begin{lemma}
Let $S$ and $U$ be elements of ${\mathcal G}_{\infty}(H)$
such that $S\dotplus U$ {\rm(}the minimal closed linear subspace containing $S+U${\rm{)}}
belongs to ${\mathcal G}_{\infty}(H)$.
Then
$$\dim(S\cap U)=\dim(f(S)\cap f(U)).$$
\end{lemma}

\begin{proof}
Lemma 4 in \cite{Pankov}.
\end{proof}

If $S$ and $U$ are closed linear subspaces of finite dimension and codimension
(respectively) then we write $[S]$ and $[U]$ for the sets of all elements of
${\mathcal G}_{\infty}(H)$ containing $S$ and contained in $U$, respectively.
Using Lemma 3, we establish the following.

\begin{lemma}
For every $S\in {\mathcal G}_{k}(H)$ and $U\in {\mathcal G}^{k}(H)$
there exist $S'\in {\mathcal G}_{k}(H)$ and $U'\in {\mathcal G}^{k}(H)$
such that
$$f([S])=[S']\;\mbox{ and }\;f([U])=[U'].$$
\end{lemma}

\begin{proof}
Lemmas 5 and 6 in \cite{Pankov}.
\end{proof}

Therefore $f$ can be extended to an automorphisms of ${\mathcal G}(H)$.
By the Fundamental Theorem of Projective Geometry,
it is induced by an invertible semi-linear operator.
Lemma 2 guarantees that the associated automorphism of ${\mathbb C}$
is identical or the complex conjugate mapping.
Since our operator preserves ${\mathcal G}^{1}(H)$,
Lemma 1 gives the claim.

\section{Remark}
Suppose that $A:H\to H$ is an invertible bounded semi-linear operator
such that $A$ and $A^{*}$ give the same transformation of ${\mathcal G}_{\infty}(H)$
which will be denoted by $f$.
In other words, $A^{*}=aA$ for a certain scalar $a$.
If $A$ is linear then
$$A= A^{**}=(aA)^{*}=\bar{a}A^{*}=|a|^{2}A$$
($a=\pm 1$ in the real case)
in the real and the complex case, respectively;
thus $|a|=1$.
If $A$ is over the complex conjugate mapping
then $(aA)^{*}=aA^{*}$ and we get $a=\pm 1$.
It is easy to see that
\begin{equation}\label{eq-1}
\forall\;S,U\in {\mathcal G}_{\infty}(H)\;\;\;\;\;S\perp f(U)\;\Longleftrightarrow\; f(S)\perp U.
\end{equation}
Now we prove the following.

\begin{lemma}
If a bijective transformation $f$ of ${\mathcal G}_{\infty}(H)$
satisfies {\rm \eqref{eq-1}}
then it is order preserving.
\end{lemma}

\begin{proof}
Since a linear subspace is contained in $S\in {\mathcal G}_{\infty}(H)$
if and only it is orthogonal to $S^{\perp}$,
the condition \eqref{eq-1} guarantees that
$f$ sends the set of all elements of ${\mathcal G}_{\infty}(H)$
contained in $S$ to the set of all elements of ${\mathcal G}_{\infty}(H)$
contained in the subspace
$$f'(S):=f^{-1}(S^{\perp})^{\perp}.$$
This means that the bijection $f'$ is order preserving.
Then $f^{-1}$ is order preserving and we get the claim.
\end{proof}

Let $f$ be as in Lemma 5.
If $H$ is real then it follows from the main result of \cite{Pankov}
that $f$ is induced by an invertible bounded linear operator $A:H\to H$.
In the complex case, Theorem 1 implies that our transformation $f$ is induced by
an invertible bounded semilinear operator $A:H\to H$ if it is continuous.
In each of these cases,
\eqref{eq-1} shows that
$$A(T)=A^{-1}(T^{\perp})^{\perp}$$
for every $T\in {\mathcal G}_{\infty}(H)$;
hence $A$ and $A^{*}$ gives the same transformation of
${\mathcal G}_{\infty}(H)$.

\end{document}